\documentclass{amsart}
\usepackage{amssymb}
\usepackage{amscd}
\usepackage{verbatim}
\usepackage{curves}

\newcommand\Real{\mathbb{R}}
\newcommand\RR{\mathbb{R}}

\newcommand\Cx{\mathbb{C}}

\newcommand\Nat{\mathbb{N}}

\newcommand\BB{\mathbb{B}}

\newcommand\HH{\mathbb{H}}
\newcommand\im{\operatorname{Im}}

\newcommand\pa{\partial}

\newcommand\sphere{\mathbb{S}}

\newcommand\ff{\operatorname{ff}}

\newcommand\cL{{\mathcal L}}

\newcommand\ha{{\frac{1}{2}}}

\newcommand\cF{{\mathcal F}}

\newcommand\cH{{\mathcal H}}

\newcommand\cR{{\mathcal R}}

\newcommand\CI{{\mathcal C}^{\infty}}
\newcommand\CIc{{\mathcal C}^{\infty}_{\text{c}}}

\renewcommand{\Box}{{\square}}
\newcommand\bcon{{\mathcal A}}

\newcommand\Diff{\operatorname{Diff}}
\newcommand\Diffb{\operatorname{Diff}_{\text{b}}}

\newcommand\rhot{\tilde\rho}

\newcommand\bl{{\text b}}

\newcommand\Tb{{}^\bl T}

\newcommand\ep{\epsilon}

\newtheorem{lemma}{Lemma}[section]
\newtheorem{prop}[lemma]{Proposition}
\newtheorem{thm}[lemma]{Theorem}

\newtheorem*{thm*}{Theorem}
\newtheorem*{prop*}{Proposition}
\newtheorem*{cor*}{Corollary}
\newtheorem*{conj*}{Conjecture}
\numberwithin{equation}{section}
\theoremstyle{remark}
\newtheorem{rem}[lemma]{Remark}
\newtheorem*{rem*}{Remark}
\theoremstyle{definition}

\newtheorem*{Def*}{Definition}

\newcommand{\bh}{\mathrm{bh}}
\newcommand{\ds}{\mathrm{dS}}
\newcommand{\oM}{\overline{\mathsf{M}}}
\newcommand{\tM}{\widetilde{\mathsf{M}}}
\newcommand{\hM}{\widehat{\mathsf{M}}}
\newcommand{\tfp}{\mathrm{tf}_+}
\newcommand{\tfm}{\mathrm{tf}_-}

\begin{document}

\title[Wave symptotics on de Sitter-Schwarzschild space]
{Asymptotics of solutions of the wave equation on de Sitter-Schwarzschild
space}
\date{12 November, 2008.}
\author[Richard Melrose]{Richard Melrose}
\address{Department of Mathematics, Massachusetts Institute of Technology,
Cambridge MA 02139-4307, U.S.A.}
\email{rbm@math.mit.edu}
\author[Ant\^onio S\'a Barreto]{Ant\^onio S\'a Barreto}
\address{Department of Mathematics, Purdue University, West Lafayette,
IN 47907-1395, U.S.A.}
\email{sabarre@math.purdue.edu}
\author[Andr\'as Vasy]{Andr\'as Vasy}
\address{Department of Mathematics, Stanford University, Stanford, CA
94305-2125, U.S.A.}
\email{andras@math.stanford.edu}
\begin{abstract} Solutions to the wave equation on de Sitter-Schwarzschild
  space with smooth initial data on a Cauchy surface are shown to decay
  exponentially to a constant at temporal infinity, with corresponding
  uniform decay on the appropriately compactified space.
\end{abstract}

\thanks{The authors gratefully acknowledge financial support for this
  project, the first from the National Science Foundation under grant
  DMS-0408993, the second under grant DMS-0500788, and the third under
  grant DMS-0201092 and DMS-0801226; they are also grateful for the
  environment at the Mathematical Sciences Research Institute, Berkeley,
  where this paper was completed.}

\maketitle

\section{Introduction}
In this paper we describe the asymptotics of solutions to the wave equation
on de Sitter-Schwarzschild space. The static model for the latter is
$M=\RR_t\times X$,
$X=(r_\bh,r_\ds)_r\times\sphere^2_\omega$ with the Lorentzian metric
\begin{equation}
g=\mu\,dt^2-\mu^{-1}\,dr^2-r^2\,d\omega^2,
\label{eq:dS-Sch-metric}\end{equation}
where
\begin{equation}
\mu=1-\frac{2m}{r}-\frac{\Lambda r^2}{3}
\label{eq:mu-def}\end{equation}
with $\Lambda$ and $m$ suitable positive constants, $0<9m^2\Lambda<1,$
$r_{\bh},r_{\ds}$ the two positive roots of $\mu$ and $d\omega^2$ the
standard metric on $\sphere^2$. We also consider the compactification
of $X$ to
$$
\bar X=[r_{\bh},r_{\ds}]_r\times\sphere^2_\omega.
$$
Then $\mu$ is a defining function for $\pa\bar X$ since it
vanishes simply at $r_{\bh},r_{\ds}$, i.e.\ $2\beta=
\frac{d\mu}{dr}\neq 0$ at $r=r_{\bh},r_{\ds}$.
Moreover, in what follows we will sometimes consider
\begin{equation}
\alpha=\mu^{\frac12}
\label{eq:alpha-def}\end{equation}
as a boundary defining function for a different compactification of $X.$ This
amounts to changing the $\CI$ structure of $\bar X$ by adjoining
$\alpha$ as a smooth function. We denote the new manifold by $X_{\ha}.$

The d'Alembertian with respect to \eqref{eq:dS-Sch-metric} is
\begin{gather}
\square= \alpha^{-2}(D_t^2-\alpha^2r^{-2}D_r(r^2\alpha^2 D_r)-\alpha^2r^{-2}
\Delta_\omega), 
\label{dal}
\end{gather}
where $\Delta_\omega$ is the Laplacian on $\sphere^2.$ We shall consider
solutions to $\Box u=0$ on $M.$

Regarding space-time as a product, up to the conformal factor $\alpha^2$,
is in fact misleading in several ways -- in particular, solutions to the wave
equation do not have simple asymptotic behavior on this space. Starting from the
stationary description of the metric, it is natural to first compactify the
time line {\em exponentially} to an interval $[0,1]_T.$ This can be done
using a diffeomorphism $T:\RR\to(0,1)$ with derivative $T'<0.$ Set
\begin{equation}\label{eq:T-def-1}
T_+=T_{\lambda,+}=e^{-2\lambda t}\ \text{in}\ t> C,
\end{equation}
with $\lambda$ to be determined and let
\begin{equation}\label{eq:T-def-2}
T=T_+\ \text{in}\ t> C.
\end{equation}
Similarly set
\begin{equation}\label{eq:T-def-3}
T_-=T_{\lambda,-}=e^{2\lambda t},\ T=1-T_-\ \text{in}\ t<-C.
\end{equation}
Near infinity $T$ depends on the free parameter $\lambda.$ The boundary
hypersurface $T_+=0$ (i\@.e\@.~$T=0$) in
\begin{equation}\label{eq:prod-comp}
[0,1]_T\times\bar X
\end{equation}
is called here the {\em future temporal face}, $T_-=0$ the
{\em past temporal face}, while $r=r_{\bh}$ and $r=r_{\ds}$ are the
{\em black hole}, resp.\ {\em de Sitter, infinity}, or together {\em
spatial infinity}.

In fact, it turns out that we need to use different values of
$\lambda$ at the two ends, $\lambda_{\bh}$ and $\lambda_{\ds}$.
This is discussed in more detail in the next section. There are product
decompositions near these boundaries
$$
[0,1]_T\times [r_{\bh},r_{\bh}+\delta)\times\sphere^2,\
[0,1]_T\times (r_{\ds}-\delta,r_{\ds}]\times\sphere^2.
$$
If $\delta$ is so large that these overlap, the transition function
is not smooth but rather is given by taking positive powers of the defining
function of the future temporal face, so the resulting space should really
be thought of having a \emph{polyhomogeneous
conormal} (but not smooth) structure in the sense of
differentiability up to the temporal faces. In particular, there
is no globally preferred boundary defining function for the temporal face,
rather such a function is only determined up to positive powers and
multiplication by positive factors. Thus, there is no fully natural `unit' of
decay but we consider powers of $e^{-t}$, resp.\ $e^t$, in a neighborhood
of the future and past temporal faces, respectively.

It turns out that there are two resolutions of this compactified space which play
a useful role in describing asymptotics. The first arises by blowing up the
corners
$$
\{0\}\times\{r_{\bh}\}\times\sphere^2,
\ \{0\}\times\{r_{\ds}\}\times\sphere^2,
\ \{1\}\times\{r_{\bh}\}\times\sphere^2,
\ \{1\}\times\{r_{\ds}\}\times\sphere^2,
$$
where the blow-up is understood to be the standard spherical blow-up
when locally the future temporal face is defined by $T_{\lambda_{\ds},+}$,
resp.\ $T_{\lambda_{\bh},+}$ at the de Sitter and black hole ends.
The resulting space is denoted $\bar M.$ The lift of the temporal and spatial
faces retain their names, while the new front faces are called the {\em
  scattering faces}. This is closely related to to the Penrose
compactification, where however the temporal faces are compressed.

Thus, a neighborhood of the lift of
$\{0\}\times\{r_{\bh}\}\times\sphere^2$ is diffeomorphic to
\begin{equation}\label{eq:bh-tfp-model}
[0,\ep)_\rho\times[r_{\bh},r_{\bh}+\delta)\times
\sphere^2_\omega,\ \rho=\rho_{\bh,+}=T_{\lambda_{\bh},+}/\mu.
\end{equation}
Similarly, a neighborhood of $\{0\}\times\{r_{\ds}\}\times\sphere^2$
is diffeomorphic to
\begin{equation}\label{eq:ds-tfp-model}
[0,\ep)_\rho\times(r_{\ds}-\delta,r_{\ds}]\times
\sphere^2_\omega,\ \rho=\rho_{\ds,+}=T_{\lambda_{\ds},+}/\mu.
\end{equation}
If $\delta>0$ is large enough, these cover a neighborhood of the
{\em future temporal
face} $\tfp$, given by the lift of $T=0$. Thus a neighborhood
of the interior of $\tfp$, is
polyhomogeneous-diffeomorphic to an open subset of
\begin{equation}\begin{split}\label{eq:int-temporal-face-model}
&[0,\ep)_x\times(r_{\bh},r_{\ds})\times
\sphere^2_\omega,\\
\qquad
&x=\rho_{\bh,+}^{1/(2\lambda_{\bh,+})}\ \text{for}\ r\ \text{near}\ r_{\bh},
\ x=\rho_{\ds,+}^{1/(2\lambda_{\ds,+})}\ \text{for}\ r\ \text{near}\ r_{\ds},
\end{split}\end{equation}
where we let the preferred defining function (up to taking
positive multiples) of $\tfp$ be $x=e^{-t}$ in the interior of $\tfp$,
hence $x=\rho_{\bh,+}^{1/(2\lambda_{\bh,+})}$ at the black hole boundary
of $\tfp$.
This means, in particular, that a neighborhood of $\tfp$ is
polyhomogeneous diffeomorphic to
\begin{equation}\label{eq:temporal-face-model}
[0,\ep)_x\times[r_{\bh},r_{\ds}]\times
\sphere^2_\omega.
\end{equation}

If $\mu$ is replaced by $\alpha$ as the defining function of the boundary
of $X$, i\@.e\@.~$\bar X$ and $T_+$ are replaced by $\bar X_{1/2}$ and
$T_+^{1/2}$ (and analogously in the past) the resulting space is denoted
$\bar M_{1/2}.$ Thus $\bar M_{1/2}$ is the square-root blow up of $\bar M,$
where the square root of the defining function of every boundary hypersurface
has been appended to the smooth structure. Here $\tfp$ is naturally
diffeomorphic to $\bar X$ in $\bar M,$ and to $\bar X_{1/2}$ in $\bar M_{1/2}.$
Both $\bar M$ and $\bar M_{1/2}$ have polyhomogeneous conormal structures
at $\tfp$ and $\tfm$; we let the preferred defining function (up to taking
positive multiples) of $\tfp$ be $x=e^{-t}$ in the interior of $\tfp$,
hence $x=\rho^{1/(2\lambda)}$ at $\pa\tfp$.

\begin{figure}[ht]
\begin{center}
\setlength{\unitlength}{.8mm}
\begin{picture}(140,80)(0,0)
\thicklines
\curve(20,20,60,20)
\curve(60,20,60,60)
\curve(60,60,20,60)
\curve(20,60,20,20)

\put(35,70){\makebox{$\bar\RR\times\bar X$}}

\curve(20,40,60,40)
\put(32,35){\makebox{$t=0$}}

\thinlines

\curve(20,47,60,47)
\curve(20,53,60,53)

\curve(30,20,30,60)
\curve(50,20,50,60)

\thicklines

\curve(90,20,120,20)
\curve(120,20,130,30)
\curve(130,30,130,50)
\curve(130,50,120,60)
\curve(120,60,90,60)
\curve(90,60,80,50)
\curve(80,50,80,30)
\curve(80,30,90,20)

\put(100,70){\makebox{$\bar M$}}
\put(105,62){\makebox{$\mathrm{tf}_+$}}
\put(105,14){\makebox{$\mathrm{tf}_-$}}

\curve(80,40,130,40)
\put(92,35){\makebox{$t=0$}}

\put(100,60){\vector(1,0){0}}
\put(83,53){\vector(-1,-1){0}}
\put(85,52){\makebox{$\rho$}}
\put(95,57){\makebox{$\mu$}}

\thinlines

\curve(80,43,105,48,130,43)
\curve(80,47,105,55,130,47)

\curve(100,20,95,40,100,60)
\curve(110,20,115,40,110,60)

\end{picture}
\end{center}
\caption{On the left, the space-time product compactification of
de Sitter-Schwarzschild space is shown (ignoring the product
with $\sphere^2$), with the time and
space coordinate lines indicated by thin lines. On the right,
$\bar M$ is shown, with the time and
space coordinates indicated by thin lines. These are no longer valid
coordinates on $\bar M$. Valid coordinates near the top left corner
are $\rho$ and $\mu$.}
\label{fig:bar-M-compactification}
\end{figure}
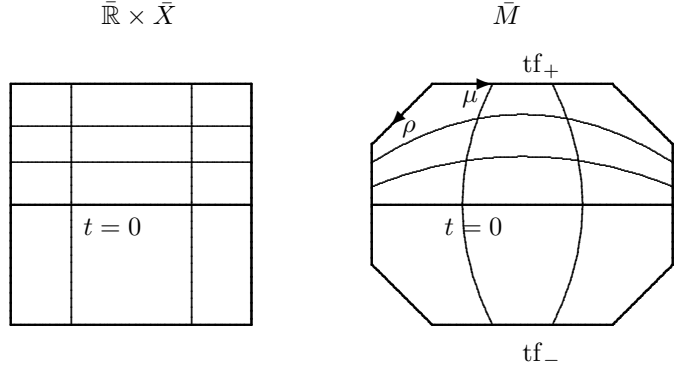

Solutions to the wave equation, when lifted to this space have simpler
asymptotics than on the product compactification,
\eqref{eq:prod-comp}. The first indication of this is that $g$ extends to
be $\CI$ and non-degenerate, up to the scattering faces, $\mu=0,$ away from
spatial infinity and uniformly up to the temporal face; the scattering
faces are characteristic with respect to the metric.  We can thus extend
$\bar M$ across $\mu=0$ to a manifold $\tilde M,$ by allowing $\mu$ to take
negative values; then the scattering face becomes an interior
characteristic hypersurface.

A further indication of the utility of this space can be seen from our main
result which is stated in terms of
\begin{equation}
\bcon^m_{\tfp}(\bar M).
\label{Con}\end{equation}
This consists of those functions which are $\CI$ on $\bar M$ away from
$\tfp,$ and  are conormal at $\tfp,$ including smoothness up to the
boundary of $\tfp.$ Such spaces are well-defined, even though the smooth
structure on $\bar M$ is not; the conormal structure suffices.
Thus, the elements of \eqref{Con}, are fixed by the condition that
for any $k$ and smooth vector fields $V_1,\ldots,V_k$ on $\bar M$ which are
tangent to $\tfp,$ 
$$
V_1\ldots V_k v\in x^m L^2_{\bl,\tfp}(\bar M),
$$
where $L^2_{\bl,\tfp}(\bar M)$ is the $L^2$-space with respect
a density $\nu_{\text{b}}$ such that $x\nu_{\text{b}}$ is smooth and strictly
positive on $\bar M$. Such a density is well-defined up to a strictly
positive polyhomogeneous multiple even under the operation of replacing $x$
by a positive power, although the weight $x^m$ is not. Thus, for all $m\in\RR,$
\begin{equation*}
x^{m+\epsilon}\CI(\bar M)\subset \bcon^{m}_{\tfp}(\bar M)
\subset x^m L^\infty(\bar M),\ \ep>0.
\end{equation*}

The main result on wave propagation is:

\begin{thm}\label{thm:asymptotics} Suppose $u\in\CI(\bar M)$ satisfies
  $\Box u=0$ for $x\in (0,1),$ then there exists a constant $c$ and
  $\ep>0$ such that
$$
u-c\in\bcon^\ep_{\tfp}(\bar M)=x^\ep  \bcon^0_{\tfp}(\bar M).
$$
\end{thm}

Thus, $u$ has an asymptotic limit, which 
happens to be a constant, at $\tfp,$ uniformly on $\bar X.$

While we have concerned ourselves with the behavior of the metric
at the corner, in regions where $\rho<C$ (i.e.\ near temporal infinity),
it is worthwhile considering what happens where $\rho>C$, i.e.\ at
spatial infinity. As we shall see, spatial infinity can be blown down,
i.e.\ there is a manifold $\oM$ and a $\CI$ map $\beta$,
$\beta:\bar M\to \oM$ such that $\beta$ is a diffeomorphism
away from spatial infinity, and such that $g$ lifts to a $\CI$ Lorentz
b-metric on $\oM,$ with tangent (i\@.e\@.~b-) behavior at the temporal
face, smooth at the other faces, 
with respect to which the non-temporal faces are characteristic.
One valid coordinate system in a neighborhood of the image of
a neighborhood of the black hole end of
spatial infinity, disjoint from temporal infinity,
is given by exponentiated versions of Eddington-Finkelstein
coordinates. In our notation, this corresponds to
$$
s_{\bh,+}=\alpha/T_{\lambda_{\bh},+}^{1/2}=\rho_{\bh,+}^{-1/2},\
 s_{\bh,-}=\alpha/T_{\lambda_{\bh},-}^{1/2}
=\alpha T_{\lambda_{\bh},+}^{1/2}=\mu \rho_{\bh,+}^{1/2},\
 \omega,
$$
where as usual $\omega$ denotes coordinates on $\sphere^2$. Here
$$
\cF_{\bh,+}=\{s_{\bh,-}=0\}
$$
is the characteristic surface given by $\mu=0$ in $T>0$ (i.e.\ the
front face of the blow up of the corner), and
$$
\cF_{\bh,-}=\{s_{\bh,+}=0\}
$$
is its
negative time analogue. The change of coordinates
$(\rho_{\bh,+},\mu)\mapsto (s_{\bh,+},s_{\bh,-})$ is a diffeomorphism
from $(0,\infty)\times (0,\delta)$ onto its image, i.e.\ these
coordinates are indeed compatible. As we show in the next section,
the metric is $\CI$ and non-degenerate on $\oM$,
and the boundary faces $s_{\bh,+}=0$ and $s_{\bh,-}=0$ are characteristic.
We can again extend $\oM$ to $\tM$, which has only two boundary faces
(the two temporal ones) by allowing $s_{\bh,\pm}$, and analogously
$s_{\ds,\pm}$, to take on negative values.
Thus, $\oM$ has six boundary faces,
$$
\tfp,\tfm,\cF_{\bh,\pm},\cF_{\ds,\pm},
$$
called the future and past temporal faces, and the future ($+$) and
past ($-$) black hole and de Sitter scattering faces.

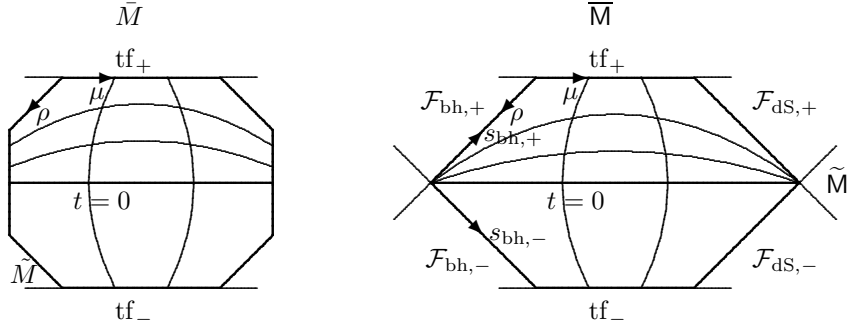
\begin{figure}[ht]
\begin{center}
\setlength{\unitlength}{0.7mm}
\begin{picture}(180,80)(0,0)

\thicklines

\curve(20,20,50,20)
\curve(50,20,60,30)
\curve(60,30,60,50)
\curve(60,50,50,60)
\curve(50,60,20,60)
\curve(20,60,10,50)
\curve(10,50,10,30)
\curve(10,30,20,20)

\put(30,70){\makebox{$\bar M$}}
\put(30,62){\makebox{$\mathrm{tf}_+$}}
\put(30,14){\makebox{$\mathrm{tf}_-$}}

\curve(10,40,60,40)
\put(22,35){\makebox{$t=0$}}

\put(30,60){\vector(1,0){0}}
\put(13,53){\vector(-1,-1){0}}
\put(15,52){\makebox{$\rho$}}
\put(25,56){\makebox{$\mu$}}

\thinlines

\curve(10,43,35,48,60,43)
\curve(10,47,35,55,60,47)

\curve(30,20,25,40,30,60)
\curve(40,20,45,40,40,60)

\curve(20,20,13,20)
\curve(50,20,57,20)
\curve(20,60,13,60)
\curve(50,60,57,60)

\put(10,21){\makebox{$\tilde M$}}

\thicklines

\curve(110,20,140,20)
\curve(140,20,160,40)
\curve(160,40,140,60)
\curve(140,60,110,60)
\curve(110,60,90,40)
\curve(90,40,110,20)

\put(120,70){\makebox{$\overline{\mathsf{M}}$}}
\put(151,54){\makebox{$\cF_{\mathrm{dS},+}$}}
\put(151,24){\makebox{$\cF_{\mathrm{dS},-}$}}
\put(88,54){\makebox{$\cF_{\mathrm{bh},+}$}}
\put(89,24){\makebox{$\cF_{\mathrm{bh},-}$}}
\put(120,62){\makebox{$\mathrm{tf}_+$}}
\put(120,14){\makebox{$\mathrm{tf}_-$}}

\put(100,50){\vector(1,1){0}}
\put(100,30){\vector(1,-1){0}}
\put(100,48){\makebox{$s_{\mathrm{bh},+}$}}
\put(101,29){\makebox{$s_{\mathrm{bh},-}$}}

\put(120,60){\vector(1,0){0}}
\put(103,53){\vector(-1,-1){0}}
\put(105,52){\makebox{$\rho$}}
\put(115,56){\makebox{$\mu$}}

\curve(90,40,160,40)
\put(112,35){\makebox{$t=0$}}

\thinlines
\curve(90,40,125,46,160,40)
\curve(90,40,125,53,160,40)

\curve(120,20,115,40,120,60)
\curve(130,20,135,40,130,60)

\curve(160,40,167,47)
\curve(160,40,167,33)
\curve(90,40,83,47)
\curve(90,40,83,33)

\curve(110,20,103,20)
\curve(140,20,147,20)
\curve(110,60,103,60)
\curve(140,60,147,60)

\put(165,38){\makebox{$\widetilde{\mathsf{M}}$}}

\end{picture}
\end{center}
\caption{On the left, $\bar M$ is shown, while on the right
its blow-down $\oM$.
The time and space coordinate lines corresponding to the product
decomposition are indicated by thin lines in the interior. The temporal
boundary hypersurfaces of $\bar M$ are continued by thin lines, as are
the characteristic surfaces $\cF_{\bh,\pm}$ and $\cF_{\ds,\pm}$, to show
that the Lorentz metric extends smoothly across $\cF_{\bh,\pm}$ and
$\cF_{\ds,\pm}$ (but not across
the temporal face!). The extended spaces are denoted by $\tilde M$ and
$\tM$.
Valid coordinates near $\cF_{\bh,+}\cap\cF_{\bh,-}$ are (apart from
the spherical coordinates) $s_{\bh,+}$ and $s_{\bh,-}$, as shown.}
\label{fig:oM-compactification}
\end{figure}

The following propagation result follows directly from the properties of
this blow-down.

\begin{prop}\label{prop:interior-propagation}
If $u$ satisfies $\Box u=0$ and has $\CI$ Cauchy data on a space-like
Cauchy surface $\Sigma\subset\tM\cap \{t\geq 0\},$ for example $\Sigma=\{t=0\}$
(i.e.\ $s_{\bh,+}=s_{\bh,-}$), then $u\in\CI(\oM^\circ).$
\end{prop}

Combining Proposition~\ref{prop:interior-propagation} and
Theorem~\ref{thm:asymptotics}, leads to the main result of this paper:

\begin{thm}\label{thm:main}
If $u$ satisfies $\Box u=0$ and has $\CI$ Cauchy data on a space-like
Cauchy surface $\Sigma\subset\tM\cap \{t\geq 0\}$ then there exists a
constant $c$ and $\ep>0$ such that 
\begin{equation}
u-c\in\bcon^\ep_{\tfp}(\oM)=x^\ep \bcon^0_{\tfp}(\oM)
\label{Limatinf}\end{equation}
near the future temporal face, $\tfp.$
\end{thm}

\begin{rem}
Our methods extend further,
for example to Cauchy data at $t=0$ which are conormal at $\pa\bar X,$ of
growth $s_{\bh,+}^m$, $m>-2$, at the boundary. Standard hyperbolic
propagation gives the same behavior at $\cF_{\bh,+}$ and $\cF_{\ds,+}$
in $\rho<C,$ and then the resolvent estimates for the `spatial Laplacian'
$\Delta_X$ (described below) apply to yield the same asymptotic term but
with convergence in the appropriate conormal space, including conormality
with respect to $\cF_{\bh,+}$ and $\cF_{\ds,+}.$ 
\end{rem}

\begin{rem} Dafermos and Rodnianski \cite{Dafermos-Rodnianski:Sch-DS}
have proved, by rather different methods, a similar result with an arbitrary
logarithmic decay rate, i\@.e\@.~an analogue of 
\begin{equation*}
u-c\in(\log\rho)^{-N}\bcon^0_{\tfp}(\oM)
\end{equation*}
for every $N.$ In terms of our approach, such logarithmic convergence 
follows from polynomial bounds on the resolvent of $\Delta_X$
{\em at the real axis}, rather than in a strip
for the analytic continuation; such estimates are much easier to obtain, as
is explained below.
\end{rem}

As already indicated, by looking at the appropriate compactification,
one only needs to study the asymptotics near $\tfp$ in $\bar M$ (or
equivalently, $\oM$). We do this by taking the Mellin transform of the wave
equation and using high-energy resolvent estimates for a `Laplacian'
$\Delta_X$ on $\bar X$. A conjugated version of this operator
is asymptotically hyperbolic, hence fits into the framework of
Mazzeo and the first author
\cite{Mazzeo-Melrose:Meromorphic}, which in particular shows
the existence of an analytic continuation for the resolvent
$$
R(\sigma)=(\Delta_X-\sigma^2)^{-1},\ \im\sigma<0.
$$
Here we also need high-energy estimates for $R(\sigma).$

The operator $\Delta_X$ has been studied by the second author
and Zworski in \cite{Sa-Barreto-Zworski:Distribution}, where it is shown
(using the spherical symmetry to reduce to a one-dimensional problem and
applying complex scaling) that the resolvent
admits an analytic continuation, from the `physical half plane',
with only one pole, at $0,$ in $\im\sigma<\ep,$ for $\ep$ sufficiently
small. Bony and Hafner in \cite{Bony-Haefner:Decay} extend and refine
this result to derive polynomial bounds on the cutoff 
resolvent, $\chi R(\sigma)\chi$, $\chi\in\CIc(X)$, as $|\sigma|\to\infty$
in the strip $|\im\sigma|<\ep.$ This implies that, for initial data in
$\CIc(X),$ the {\em local energy}, i.e.\ the energy in a fixed compact set
in space, decays to the energy corresponding to the 0-resonance. In our
terminology this amounts to studying the behavior 
of the solution near a compact subset of the interior of $\tfp.$ Our
extension of their result is both to allow more general initial data, not
necessarily of compact support, and to study the asymptotics uniformly up to
the boundary at temporal infinity. This requires resolvent estimates on slightly
weighted $L^2$-spaces, which were obtained by the authors in
\cite{Melrose-SaBarreto-Vasy:Semiclassical} together with the use of the
geometric compactification $\bar M$ (or $\oM).$ For this to succeed, it is
essential that the resolvent only be applied to `errors' which intersect
$\pa\oM$ in the interior of $\cF_{\bh,+}$ and $\cF_{\ds,+}$. This turns out
to be a major gain since the analytic continuation of the resolvent (even
arbitrarily close to the real axis) cannot be applied directly to the
initial data. Thus essential use is made of the fact that once the solution
has been propagated to the scattering faces, the error terms have more decay.

It is then relatively clear, as remarked above, that if one only knew
polynomial growth estimates for the limiting resolvent {\em at the real
  axis} (rather than in a strip), one could still obtain the same
asymptotics, but with error that is only super-logarithmically decaying.
This observation may be of use in other settings where such polynomial
bounds are relatively easy to obtain from estimates for the cutoff
resolvent, as in \cite{Bony-Haefner:Decay}, or analogous semiclassical
propagation estimates at the trapped set, by pasting with well-known high
energy resolvent estimates localized near infinity. This has been studied
particularly by Cardoso and Vodev \cite{Cardoso-Vodev:Uniform}, using the
method of Bruneau and Petkov
\cite[Section~3]{Bruneau-Petkov:Semiclassical}.

This paper is structured as follows. In Section~\ref{section:geometry} both
the compactifications and the underlying geometry are discussed in more
detail. The `spatial Laplacian' and relevant resolvent estimates are
recalled in Section~\ref{section:resolvent} and in
Section~\ref{section:asymptotics} the main result is proved using the
Mellin transform.

\section{Geometry}\label{section:geometry} In this section the
various compactifications of de Sitter-Schwarzschild space are studied
after an initial examination of the simpler case of de Sitter space.

\subsection{De Sitter space}
We start with the extreme case of de Sitter space, corresponding
to $m=0$ in \eqref{eq:dS-Sch-metric} and \eqref{eq:mu-def}, to see what the
`correct' compactification of $M$ should be. However, rather than starting
from the static model, consider this as a Lorentzian symmetric space.
De Sitter space is given by the hyperboloid
$$
z_1^2+\ldots+z_n^2=z_{n+1}^2+1\ \text{in}\ \Real^{n+1}
$$
equipped with the pull-back of the Minkowski metric
$$
dz_{n+1}^2-dz_1^2-\ldots-dz_n^2.
$$
Introducing polar coordinates $(R,\theta)$ in $(z_1,\ldots,z_n)$,
so
$$
R=\sqrt{z_1^2+\ldots+z_n^2}=\sqrt{1+z_{n+1}^2},
\ \theta=R^{-1}(z_1,\ldots,z_n)\in\sphere^{n-1}
,\ \tau=z_{n+1},
$$
the hyperboloid
can be identified with $\Real_\tau\times\sphere^{n-1}_\theta$ with the
Lorentzian metric
\begin{equation*}
\frac{d\tau^2}{\tau^2+1}-(\tau^2+1)\,d\theta^2,
\end{equation*}
where $d\theta^2$ is the standard Riemannian metric on the sphere.
For $\tau>1,$ set $x=\tau^{-1}$, so the metric becomes
$$
\frac{(1+x^2)^{-1}\,dx^2-(1+x^2)\,d\theta^2}{x^2}.
$$
An analogous formula holds for $\tau<-1$, so compactifying
the real line to an interval $[0,1]_T$, with $T=x=\tau^{-1}$
for $x<\frac{1}{4}$ (i.e.\ $\tau>4$), say, and $T=1-|\tau|^{-1}$, $\tau<-4$,
gives a compactification, $\hM,$ of
de Sitter space on which the metric is conformal to a non-degenerate
Lorentz metric. There is natural generalization, to 
{\em asymptotically de Sitter-like spaces} $\hM$,
which are diffeomorphic to compactifications $[0,1]_T\times Y$
of $\RR_\tau\times Y$, where $Y$ is a compact manifold
without boundary, and $\hM$ is equipped
with a Lorentz metric on its interior which is conformal
to a Lorentz metric smooth up to the boundary. These
space-times are Lorentzian analogues of the much-studied conformally
compact (Riemannian) spaces. On this class of space-times the solutions of
the Klein-Gordon equation were analyzed by the third author
in \cite{Vasy:De-Sitter}, and were shown to have simple asymptotics
analogous to those for eigenfunctions on conformally compact manifolds.

\begin{thm*}(\cite[Theorem~1.1.]{Vasy:De-Sitter})
Set $s_\pm(\lambda)=\frac{n-1}{2}\pm\sqrt{\frac{(n-1)^2}{4}-\lambda}.$
If $s_+(\lambda)-s_-(\lambda)\notin\Nat,$ any solution $u$ of the Cauchy
problem for $\Box-\lambda$ with $\CI$ initial data at $\tau=0$ is of the form
\begin{equation*}
u=x^{s_+(\lambda)}v_++x^{s_-(\lambda)}v_-,\ v_\pm\in\CI(\hM).
\end{equation*}
If $s_+(\lambda)-s_-(\lambda)$ is an integer, the same conclusion holds if
$v_-\in\CI(\hM)$ is replaced by $v_-=\CI(\hM)
+x^{s_+(\lambda)-s_-(\lambda)}\log x\,\CI(\hM)$.
\end{thm*}

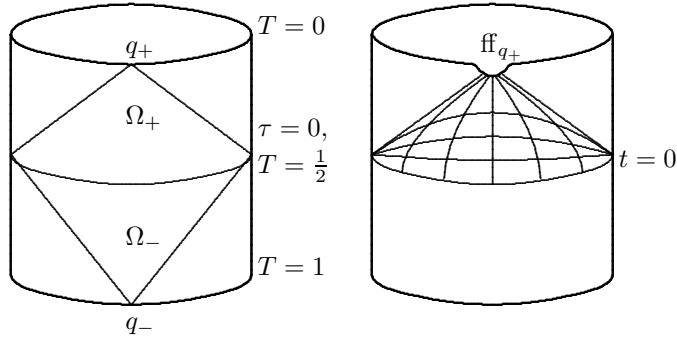
\begin{figure}[ht]
\begin{center}
\setlength{\unitlength}{0.8mm}
\begin{picture}(140,75)(0,0)

\thicklines

\curve(20,20,25,17,40,15,55,17,60,20)
\curve(20,60,25,57,40,55,55,57,60,60)
\curve(20,60,25,63,40,65,55,63,60,60)
\curve(20,20,20,60)
\curve(60,20,60,60)

\thinlines
\curve(20,40,25,37,40,35,55,37,60,40)

\curve(40,55,20,40)
\curve(40,55,60,40)
\curve(40,15,20,40)
\curve(40,15,60,40)

\put(61,20){\makebox{$T=1$}}
\put(61,43){\makebox{$\tau=0,$}}
\put(61,37){\makebox{$T=\frac12$}}
\put(61,60){\makebox{$T=0$}}
\put(39,57){\makebox{$q_+$}}
\put(39,11){\makebox{$q_-$}}
\put(39,45){\makebox{$\Omega_+$}}
\put(39,25){\makebox{$\Omega_-$}}


\thicklines

\curve(80,20,85,17,100,15,115,17,120,20)
\curve(80,60,85,57,97,55)
\curve(97,55,100,53,103,55)
\curve(103,55,115,57,120,60)
\curve(80,60,85,63,100,65,115,63,120,60)
\curve(80,20,80,60)
\curve(120,20,120,60)

\thinlines
\curve(80,40,85,37,100,35,115,37,120,40)

\put(121,38){\makebox{$t=0$}}

\curve(98.3,53.3,80,40)
\curve(101.7,53.3,120,40)

\put(98,57){\makebox{$\mathrm{ff}_{q_+}$}}

\curve(100,53,100,35)
\curve(101,53.1,111.5,44,115,37)
\curve(99,53.1,88.5,44,85,37)
\curve(99.5,53,94,44,92,36)
\curve(100.5,53,106,44,108,36)

\curve(80,40,100,39,120,40)
\curve(80,40,100,43,120,40)
\curve(80,40,100,47,120,40)

\end{picture}
\end{center}
\caption{On the left, the compactification of de Sitter space with the
backward light cone from $q_+=(1,0,0,0)$ and forward light cone from $q_-
=(-1,0,0,0)$ are
shown. $\Omega_+$, resp.\ $\Omega_-$, denotes the intersection of these
light cones with $t>0$, resp.\ $t<0$.
On the right, the blow up of de Sitter space at $q_+$ is shown. The
interior of the light cone inside the front face $\ff_{q_+}$ can
be identified with the spatial part of the static model of de Sitter
space. The spatial and temporal coordinate lines for the static model
are also shown.}
\label{fig:deSitter-space}
\end{figure}

The simple structure of the de Sitter metric (and to some extent
the asymptotically de Sitter-like metrics) can be hidden by
blowing up certain submanifolds of $\hM$.
In particular, the {\em static model} of de Sitter space arises by
singling out a point on
$\sphere^{n-1}_{\theta}$,
e.g.\ $q_0=(1,0,\ldots,0)\in\sphere^{n-1}\subset\RR^n$.
Note that $(\theta_2,\ldots,\theta_n)\in\RR^{n-1}$ are local
coordinates on $\sphere^{n-1}$ near $q_0$.
Now consider
the intersection of the backward light cone
from $q_0$ considered as a point $q_+$ at future infinity, i.e.\ 
where $T=0$, and the forward
light cone from $q_0$ considered as a point $q_-$ at past infinity,
i.e.\ where $T=1$.
These intersect the equator $T=1/2$ (here $\tau=0$)
in the same set, and together form a `diamond', $\hat\Omega$, with a conic
singularity at $q_+$ and $q_-.$ Explicitly $\hat\Omega$
is given by $z_2^2+\ldots+z_n^2\leq 1$ inside the
hyperboloid. If $q_+$, $q_-$ are blown up, as well as
the corner $\pa\Omega\cap\{\tau=0\}$, i.e.\ where
the light cones intersect $\tau=0$ in $\hat\Omega$,
we obtain a manifold
$\bar M$,
which can be blown down to (i.e.\ is a blow up of) the space-time
product $[0,1]\times\overline{\BB^{n-1}}$, with
$\BB^{n-1}=\{Z\in\RR^{n-1}:\ |Z|<1\}$ on which the Lorentz metric has
a time-translation invariant warped product form. Namely, first considering
the interior $\Omega$ of $\hat\Omega$ we introduce
the global (in $\Omega$) standard static
coordinates $(t,Z)$, given by (with the expressions
involving $x$ valid near $T=0$)
\begin{equation*}\begin{split}
\BB^{n-1}\ni Z&=(z_2,\ldots,z_n)=x^{-1}\sqrt{1+x^2}(\theta_2,\ldots,\theta_n),\\
\sinh t&=\frac{z_{n+1}}{\sqrt{z_1^2-z_{n+1}^2}}=(x^2-(1+x^2)(\theta_2^2+
\ldots+\theta_n^2))^{-1/2},
\end{split}\end{equation*}
It is convenient to rewrite these as well in terms of polar coordinates
in $Z$ (valid away from $Z=0$):
\begin{equation*}\begin{split}
r&=\sqrt{z_2^2+\ldots+z_n^2}=\sqrt{1+z_{n+1}^2-z_1^2}
=x^{-1}\sqrt{1+x^2}\sqrt{\theta_2^2+\ldots
+\theta_n^2},\\
\sinh t&=\frac{z_{n+1}}{\sqrt{z_1^2-z_{n+1}^2}}=(x^2-(1+x^2)(\theta_2^2+
\ldots+\theta_n^2))^{-1/2}=x^{-1}(1-r^2)^{-1/2},\\
&\ \omega=r^{-1}(z_2,\ldots,z_n)=(\theta_2^2+\ldots+\theta_n^2)^{-1/2}
(\theta_2,\ldots,\theta_n)\in\sphere^{n-2}.
\end{split}\end{equation*}
In these coordinates the
metric becomes
$$
(1-r^2)\,dt^2-(1-r^2)^{-1}dr^2-r^2\,d\omega^2,
$$
which is a special case of the de Sitter-Schwarzschild metrics with $m=0$
and $\Lambda=3$.

\begin{lemma}
The lift of $\hat\Omega$ to the blow up
$[\hM;q_+,q_-]$ is a $\CI$ manifold
with corners, $\bar\Omega$. Moreover,
$[\bar\Omega,\pa\Omega\cap\{\tau=0\}]$ is naturally
diffeomorphic to the $\CI$ manifold with corners obtained from
$[[0,1]\times\overline{\BB^{n-1}};\{0\}\times\pa\BB^{n-1};
\{1\}\times\pa\BB^{n-1}]$ by adding the square root of the defining function
of the lift of $\{0\}\times\overline{\BB^{n-1}}$ and $\{1\}\times\overline{\BB^{n-1}}$
to the $\CI$ structure.
\end{lemma}

\begin{rem}
This lemma states that from the stationary
point of view, the `right' compactification near the top face arises
by blowing up the corner $\pa\overline{\RR}\times\pa\overline{\BB^{n-1}}$,
although the resulting space is actually more complicated than needed,
since the original boundary hypersurfaces of the stationary space
can be blown down to obtain a subset of $\hM$.

The fact that this approach gives $\rho=x^2$ as the
defining function of the temporal face, rather than $x$ (hence necessitating
adding the square root of $\rho$ to the smooth structure),
corresponds to the fact that,
in the sense of Guillarmou \cite{Guillarmou:Meromorphic},
$\hM$ is actually {\em even}.
\end{rem}

\begin{proof}
Coordinates on the blow up of $\hM$ near the lift of $q_+$ are given by
$$
x,\ \Theta_j=\frac{\sqrt{1+x^2}}{x}\,\theta_j,\ j=2,\ldots,n;
$$
these are all bounded in the region of validity of the coordinates, with
the light cone given by $\sum_{j=2}^n\Theta_j^2=1$ (which is
why the factor $\sqrt{1+x^2}$ was introduced), so the lift $\bar\Omega$
of $\hat\Omega$
is $\sum_{j=2}^n\Theta_j^2\leq 1$. As $\sum_{j=2}^n\Theta_j^2-1$ has a
non-vanishing differential where it vanishes, this shows that
$\bar\Omega$ is a $\CI$ manifold
with corners. Near the light cone, $\pa\bar\Omega$ one
can introduce polar coordinates in $\Theta$, and use
$$
x,\ r=(\sum_{j=2}^N\Theta_j^2)^{1/2},\ (\omega_2,\ldots,\omega_n)=
r^{-1}(\Theta_2,\ldots,\Theta_n)\in\sphere^{n-1}
$$
as local coordinates.
On the other hand, blowing up the corner of $[0,1]_T\times\overline{\BB^{n-1}}$, where
$T=e^{-2t}$ for $t>4$, say, which is equivalent to $(\sinh t)^{-2}$ there,
gives coordinates near the lift of $T=0:$
$$
r,\ \rho= (\sinh t)^{-2}/(1-r^2)=x^2,\ \omega.
$$
Thus, one {\em almost} has a diffeomorphism between the two coordinate
charts, hence locally between the manifolds, except that in the
blow up of $[0,1]_T\times\overline{\BB^{n-1}}$,
the defining function of the temporal
face is the square of the defining function of the temporal face
arising from the blow up of $\hM$. This is remedied by
adding the square root of the defining function of the lift of
$\{0\}\times\overline{\BB^{n-1}}$ and $\{1\}\times\overline{\BB^{n-1}}$, i.e.\ of
$x^2$, to the smooth structure, thus proving the lemma.
\end{proof}

It is worthwhile comparing the de Sitter space wave asymptotics,
which is
\begin{equation*}
u=x^{n-1}v_++v_-,\ v_+\in\CI(\hM),\ v_-\in\CI(\hM)+x^{n-1}(\log x)\CI(\hM),
\end{equation*}
with our main result. The fact that the coefficients in the de Sitter expansion
are $\CI$ on $\hM$ means that on $\bar M$, the leading terms are
constant. Thus, the de Sitter result implies (and is much stronger than)
the statement that $u$ decays to a constant on $\bar M$ at an exponential
rate.

\subsection{Blow-up of the space-time product} We now turn to the
compactification of de Sitter-Schwarzschild space. It turns out that while
this cannot be embedded into a space as simple as $\hM$ the final setting
is not much more complicated. In terms of de Sitter space, the difference
is that while spatial infinity in $[0,1]_T\times\overline{\BB^3},$ blown up
at the corner, can always be blown down, the same is not true for temporal
infinity.

In fact, the `black hole end' $r=r_{\bh}$ resembles the `de Sitter end'
quite closely, which motivates the construction in the de
Sitter-Schwarzschild setting. There is a simpler construction, depending
on the choice of a constant $\lambda>0$, which does
not quite work because of some incompatibility between the two ends (which
whilst very similar, are not quantitatively the same). With $\mu$ as in
\eqref{eq:mu-def}, compactify $M$ by compactifying $\RR$ into $[0,1]_T$ as
in \eqref{eq:T-def-1}-\eqref{eq:T-def-3}, so
\begin{equation*}
T_+=T_{\lambda,+}=e^{-2\lambda t}\ \text{in}\ t>c,
\end{equation*}
and compactify $(r_{\bh},r_{\ds})$ as
$[r_{\bh},r_{\ds}]$, to obtain
$$
[0,1]_{T}\times [r_{\bh},r_{\ds}]\times\sphere^2=[0,1]_T\times\bar X.
$$
Then blowing up the corners
$$
\{0\}\times\{r_{\bh}\}\times\sphere^2,\ \{0\}\times\{r_{\ds}\}\times\sphere^2
$$
(and analogously at $T=1,$ i.e.\ $T_-=0$), gives a space denoted $\bar M.$ 
Thus, a neighborhood $U=U_{\lambda,+}$ of the `temporal
face' $T_+=0$ is diffeomorphic to
\begin{equation}\label{eq:temporal-face-model-2}
[0,\ep)_{\rho}\times[r_{\bh},r_{\ds}]\times
\sphere^2_\omega,\ \rho=T_{\lambda,+}/\mu.
\end{equation}
In the interior of the temporal face, where $\mu>0$, this is
in turn diffeomorphic to an open subset of
\begin{equation}
[0,\ep)_T\times(r_{\bh},r_{\ds})\times
\sphere^2_\omega.
\end{equation}

If the same construction is performed, but using the smooth structure
on the the compactification of $X$ given by $\alpha=\mu^{1/2},$
i.e.\ $\bar X_{1/2},$ and
$$
\tilde T_{\lambda,+}=T_{\lambda,+}^{1/2}=e^{-\lambda t},\ t>c,
$$
then a neighborhood $U$ of $\tfp$ as above is diffeomorphic to
\begin{equation}
[0,\ep)_{\rhot}\times[r_{\bh},r_{\ds}]_{1/2}\times
\sphere^2_\omega,\ \rhot=\tilde T_{\lambda,+}/\alpha,
\label{eq:temporal-face-model-3}\end{equation}
where $[r_{\bh},r_{\ds}]_{1/2}$ denotes that $\alpha$ has been added
to the smooth structure (or equivalently
$(r-r_{\bh})^{1/2}$ and $(r_{\ds}-r)^{1/2}$ have been added to the smooth
structure).
The distinction between \eqref{eq:temporal-face-model-2}
and \eqref{eq:temporal-face-model-3} is the same as between
$[0,1]_T\times \bar X$ and $[0,1]_{\tilde T}\times\bar X_{1/2}$
(where $\tilde T$ is defined analogously to $T$), namely
the square roots of the defining functions of all boundary hypersurfaces
have been added to the smooth structure. We denote the resulting
space by $\bar M_{1/2}$.

The subtlety is that the de Sitter and black hole ends need different
values of $\lambda.$ So what we actually need is to paste together
$U_{\lambda_{\bh},+}\cap\{r<r_2\}$ and $U_{\lambda_{\ds},+}\cap\{r>r_1\}$ for
some $r_1,r_2$, $r_{\bh}<r_1<r_2<r_{\ds}$ where
$\lambda_{\bh}$ and $\lambda_{\ds}$
are chosen in a way that reflects the local geometry neat the two ends.
The transition map in the overlap, where $r\in(r_1,r_2)$,
is given by
$$
(\rho_{\bh,+},r,\omega)\mapsto(\rho_{\ds,+},r,\omega),
\ \rho_{\ds,+}=\mu^{\lambda_{\ds}/\lambda_{\bh}-1}
\,\rho_{\bh,+}^{\lambda_{\ds}/\lambda_{\bh}}.
$$
In the overlap, where $\mu\neq 0$ so
$\mu\mapsto \mu^{\lambda_{\ds}/\lambda_{\bh}-1}$ is smooth,
this is a polyhomogeneous-diffeomorphism, in the sense that
it is polyhomogeneous in the local defining functions of the lift of
$T=0$ (namely $\rho_{\bh,+}$ and $\rho_{\ds,+}$). In particular,
the front faces of the blow-ups have well defined boundary defining functions,
namely $\mu$, up to multiplying by a $\CI$ non-vanishing function, so
we consider the resulting space a polyhomogeneous manifold with corners,
where the `polyhomogeneous' (as opposed to $\CI$) faces are only
the future and past temporal faces, $\tfp$ and $\tfm$. There is also
an analogous construction for $\bar M_{1/2}$.

As indicated already, we also want a preferred defining function (up to taking
positive multiples) $x$ of $\tfp$ in order to measure the rate of
decay at the temporal faces; this should be polyhomogoneous-equivalent
to the local defining functions $\rho_{\bh,+}$ and $\rho_{\ds,+}$.
We take this to
be of the form $x=f(r)e^{-t}$, $f>0$ smooth
for $r\in(r_{\bh},r_{\ds})$. Comparison with
\eqref{eq:bh-tfp-model}-\eqref{eq:ds-tfp-model} shows that we need to take
$f(r)=\mu^{-1/(2\lambda_{\bh})}$ for $r$ near $r_{\bh}$,
$f(r)=\mu^{-1/(2\lambda_{\ds})}$ for $r$ near $r_{\ds}$,
hence $x=\rho_{\bh,+}^{1/(2\lambda_{\bh})}$ near the black hole boundary
of $\tfp$, and $x=\rho_{\bh,+}^{1/(2\lambda_{\ds})}$ near the de Sitter
boundary.
Then, in particular, a neighborhood of $\tfp$ is
polyhomogeneous diffeomorphic to
\begin{equation}
[0,\ep)_x\times[r_{\bh},r_{\ds}]\times
\sphere^2_\omega.
\end{equation}

We still need to determine the values $\lambda$ at the two ends.
Writing $T=T_+=T_{\lambda,+}$, the dual metric (which is the principal
symbol of the wave operator, $\Box$)
has the form
\begin{equation*}
G=4\alpha^{-2}\lambda^2 T^2\pa_T^2-\alpha^2\pa_r^2-r^{-2}\pa_{\omega}^2
\end{equation*}
in the original product compactification, with
$\pa_r=\frac{d\mu}{dr}\pa_\mu=2\beta\pa_\mu$.
The change of variables from $r$ to $\mu$ is smooth and non-degenerate,
i.e.\ $2\beta=d\mu/dr\neq 0$ for $\mu$ close to $0$, i.e.\ $r$ close to
$r_{\bh}$ or $r_{\ds}$.
Note that $\beta>0$ for $r$ near $r_{\bh}$, $\beta<0$ for
$r$ near $r_{\ds}$ since $\mu>0$ for $r\in(r_{\bh},r_{\ds})$.

After blow-up, in the coordinates $(\rho,\mu,\omega),$
\begin{equation*}
G=4\mu^{-1}\lambda^2 \rho^2\pa_\rho^2
-4\mu \beta^2(\pa_\mu-\mu^{-1}\rho\pa_\rho)^2
-r^{-2}\pa_{\omega}^2.
\end{equation*}
Thus,
\begin{equation*}
G=\mu^{-1}\left(4\lambda^2 \rho^2\pa_\rho^2
-4\beta^2(\mu\pa_\mu-\rho\pa_\rho)^2\right)
-r^{-2}\pa_{\omega}^2.
\end{equation*}
If we set $\lambda=\beta(r_{\bh})>0$ or $\lambda=-\beta(r_{\ds})>0$
then\footnote{In terms of the `spatial Laplacian',
$\Delta_X$, described in Section~\ref{section:resolvent}, $-\beta(r_{\bh})^2$
and $-\beta(r_{\ds})^2$ are the asymptotic curvatures, hence they
are natural quantities from that point of view as well.}
the $\rho^2\pa_\rho^2$ terms cancel, so {\em locally near $r_{\bh}$}
\begin{equation*}
G=4\gamma\rho^2\pa_\rho^2+8\beta^2\rho\pa_\rho\pa_\mu
-4\beta^2\mu\pa_\mu^2-r^{-2}\pa_{\omega}^2,\ \gamma=\mu^{-1}(\beta(r_{\bh})^2
-\beta^2),
\end{equation*}
where $\gamma$ is $\CI$ by Taylor's theorem, and there is
a similar expansion at $r_{\ds}.$ Thus the choice of $\lambda$ determines
the compactification $\bar M,$ and it is only at this point that the
compactification has been specified. Note that this metric is a $\CI$ Lorentzian
b-metric on $[0,\ep)_\rho\times (r_{\bh}-\ep,r_{\bh}+\ep)_r\times\sphere^2_\omega$
(i.e.\ is non-degenerate as a quadratic form on the b-cotangent bundle),
in particular it is $\CI$ across $\mu=0$. Denoting
this extension of $\bar M$ by $\tilde M$ (which is now non-compact); $g$ becomes
a polyhomogeneous conormal Lorentz metric on $\tilde M$, smooth near
$\mu=0$ (where there is a well-defined smooth structure). We write $\cF$
for the set given by $\mu=0,$ i.e.\ the boundary hypersurface of $\bar M$
that is no longer a boundary hypersurface of $\tilde M$.

For this metric $\cF$ is characteristic, and one has the standard
propagation of singularities in $\rho>0$. In particular, for $\CI$ initial
data the solution is smooth in $\rho>0$ across $\mu=0.$ In fact, writing
covectors as $\xi\,\frac{d\rho}{\rho} +\zeta\,d\mu+\sum\eta_j\,d\omega_j$,
i.e.\ $(\rho,\mu,\omega,\xi,\zeta,\eta)$ are coordinates on $\Tb^*\tilde
M$, the dual metric, considered as a function on $\Tb^*\tilde M$, is
\begin{equation*}
G=4\gamma\xi^2+8\beta^2\xi\zeta
-4\beta^2\mu\zeta^2-r^{-2}|\eta|^2,
\end{equation*}
so the Hamilton vector field of $G$ is
\begin{equation*}\begin{split}
H_G=&8(\gamma\xi+\beta^2\zeta)\rho\pa_\rho-8\beta^2(\mu\zeta-\xi)\pa_\mu\\
&
-\big(4\frac{\pa\gamma}{\pa\mu}\xi^2+8\beta\frac{\pa\beta}{\pa\mu}(2\xi\zeta
-\mu\zeta^2)-4\beta^2\zeta^2-\frac{\pa r^{-2}}{\pa\mu}|\eta|^2\big)\pa_\zeta
-r^{-2}H_{(\omega,\eta)},
\end{split}\end{equation*}
with $H_{(\omega,\eta)}$ denoting the Hamilton vector field of the standard
metric on the sphere.
The conormal bundle $N^*\{\mu=0\}$ is $\mu=0,$ $\xi=0,$
$\eta=0,$ so at this set
\begin{equation*}
H_G=8\beta^2\zeta\rho\pa_\rho+4\beta^2\zeta^2\pa_\zeta
\end{equation*}
and so is indeed tangent to $N^*\{\mu=0\},$ and it is non-radial off the
zero section (where $\zeta\neq 0$) as long as $\rho\neq 0$.

At $\pa \cF$, i.e.\ at $\rho=0$, however there are radial points over the
conormal bundle of $\cF$. Rather than dealing with them directly, which can
be done in the spirit of \cite{Hassell-Melrose-Vasy:Microlocal}, we reduce
the problem to the study of the high energy behavior of the resolvent of
the spatial Laplacian (which gives more, in fact), which was performed in
\cite{Melrose-SaBarreto-Vasy:Semiclassical}.

\subsection{Blow down of spatial infinity}

We now discuss the manifold with corners $\oM$,
in which spatial infinity is blown down.
A valid coordinate system near the image of the black hole end of
spatial infinity, disjoint from temporal infinity, is given by
$$
s_{\bh,+}=\alpha/T_{\lambda_{\bh},+}^{1/2}=\rho_{\bh,+}^{-1/2},
\ s_{\bh,-}=\alpha/T_{\lambda_{\bh},-}^{1/2}
=\alpha T_{\lambda_{\bh},+}^{1/2}=\mu \rho_{\bh,+}^{1/2},
\ \omega,
$$
where as usual $\omega$ denotes coordinates on $\sphere^2$.
In these coordinates the dual metric is
$$
G=\gamma(s_{\bh,+}^2\pa_{s_{\bh,+}}^2+s_{\bh,-}^2\pa_{s_{\bh,-}}^2)
-2(\beta(r_{\bh})^2+\beta^2)\pa_{s_{\bh,+}}\pa_{s_{\bh,-}}
-r^{-2}\pa_\omega^2,
$$
and the boundary faces $s_{\bh,+}=0$ and $s_{\bh,-}=0$ are characteristic.
We can also extend $G$ to a smooth non-degenerate Lorentz metric on
$\tM^\circ$. Recall that $\tM$ contains $\oM$ as a closed domain
with corners, namely locally in $\tM$ we simply allow the four boundary
defining functions $s_{\bh,\pm}$, $s_{\ds,\pm}$, to assume negative values
(but we do not extend $\oM$ across the temporal faces).
This calculation shows the following important fact:

\begin{lemma}\label{lemma:Box-form} The d'Alembertian
  $\Box\in\Diff^2(\tM^\circ)$, and indeed $\Box\in\Diffb^2(\tM)$. Moreover
  the scattering surfaces $\cF_{\bh,\pm}=\{s_{\bh,\mp}=0\}$ and
  $\cF_{\ds,\pm}=\{s_{\ds,\mp}=0\}$ are characteristic.
\end{lemma}

Proposition~\ref{prop:interior-propagation} is an immediate
corollary of this lemma and standard hyperbolic propagation.

\section{Resolvent estimates}\label{section:resolvent}

Next consider the `spatial Laplacian', resolvent estimates for which constitute
one of the key ingredients in our analysis. From \eqref{dal} with
$T=e^{-t},$ it follows that
\begin{equation*}
\Box=\alpha^{-2}\left((TD_T)^2-\alpha^2r^{-2}D_r\alpha^2 r^{2}D_r
-\alpha^2r^{-2}\Delta_\omega\right).
\end{equation*}
Recall that not precisely $T,$ but rather $T^{\lambda_{\bh}}$ and
$T^{\lambda_{\ds}},$ were used above to construct the compactification.

By definition the spatial `Laplacian' is
\begin{equation*}\begin{split}
\Delta_X&=\alpha^2r^{-2}D_r\alpha^2 r^{2}D_r
+\alpha^2r^{-2}\Delta_\omega
\end{split}\end{equation*}
Near $\alpha=0,$ where $\alpha$ can be used as a valid coordinate
in place of $r,$
\begin{equation*}
\Delta_X
=\beta r^{-2}\alpha D_\alpha \beta r^2\alpha D_\alpha+\alpha^2r^{-2}\Delta_\omega
\in \Diff^2_0(\bar X_{1/2}).
\end{equation*}

This is not the Laplacian of a Riemannian metric on $X;$ however it is
very similar to one. It is of the form $d^*d$ with respect to the inner product
on one-forms given by the fiber inner product with respect to the
`spatial part' $H=\alpha^2\pa_r^2+r^{-2}\pa^2_\omega$ of $G$ and density
on $X$ given by $dh=\alpha^{-2} r^2\,dr\,d\omega.$ In what follows we will
also view  $\Delta_X$ is a 0-operator on $\bar X_{1/2},$
self-adjoint on
$$
L^2(X,|dh|),\ |dh|=\alpha^{-2} r^2\,|dr|\,|d\omega|
=\alpha^{-1}|\beta|^{-1}r^2\,|d\alpha|\,|d\omega|,
$$
and we will use the techniques of \cite{Mazzeo-Melrose:Meromorphic}
to study its resolvent. It is also useful to introduce the operator
$$
L=\alpha\Delta_X\alpha^{-1},
$$
which is self-adjoint on
$$
L^2(X,\alpha^{-2}\,|dh|)=\alpha L^2(X,dh),\ \alpha^{-2}|dh|
=\alpha^{-3}|\beta|^{-1}r^2\,|d\alpha|\,|d\omega|.
$$
Thus, this space is $L^2_0(X)$ as a Banach space, up to equivalence of
norms.

Let $\tilde\alpha=\alpha^{1/\lambda_{\bh}}\in\CI(X)$, $\tilde\alpha>0,$ for
$r$ near $r_{\bh}$ or $\tilde\alpha=\alpha^{1/\lambda_{\ds}}$ for $r$ near
$r_{\ds}.$ The normal operators $N_{0,\bh}(L)$, $N_{0,\ds}(L)$ of $L$ in
$\Diff^2_0(\bar X_{1/2})$ at $r=r_{\bh}$, resp.\ $r=r_{\ds}$, are
$$
N_{0,\bh}(L)=\lambda_{\bh}^2 N_{0,\bh}(\Delta_{\HH^3}),
\ N_{0,\ds}(L)=\lambda_{\ds}^2 N_{0,\ds}(\Delta_{\HH^3}),
$$
where $\Delta_{\HH^3}$ is the hyperbolic
Laplacian, explaining the usefulness of this conjugation. In particular,
it follows immediately from \cite{Mazzeo-Melrose:Meromorphic}
(with improvements
from \cite{Guillarmou:Meromorphic})
that the resolvent
$$
\cR(\sigma)=(L-\sigma^2)^{-1},\ \text{on}\ L^2_0(\bar X_{1/2}),\ \im\sigma<0,
$$
continues meromorphically
to a strip $|\im\sigma|<\ep$ as an operator between weighted
$L^2$-spaces (as well as other spaces):
$$
\cR(\sigma):\tilde\alpha^\delta L^2_0(\bar X_{1/2})
\to \tilde\alpha^{-\delta} L^2_0(\bar X_{1/2}),\ \delta>\ep;
$$
we keep denoting the analytic continuation
by $\cR(\sigma)$.
Thus,
$$
R(\sigma)=(\Delta_X-\sigma^2)^{-1}=\alpha^{-1}\cR(\sigma)\alpha\ \text{on}
\ L^2(X,|dh|),\ \im\sigma<0,
$$
continues meromorphically to a strip $|\im\sigma|<\ep$
$$
R(\sigma):\tilde\alpha^\delta L^2(X,|dh|)
\to \tilde\alpha^{-\delta} L^2(X,|dh|),\ \delta>\ep.
$$

The result we need is proved in
\cite{Melrose-SaBarreto-Vasy:Semiclassical}, giving polynomial bounds on
the resolvent in a strip around the real axis.

\begin{prop}\label{bounds} If $\ep>0$ is sufficiently small the only pole of
  the analytic continuation of the resolvent $R(\sigma)$ in $\im\sigma<\ep$
  is $\sigma=0,$ which is simple, with residue given by a constant $\gamma$ and
  for each $k$ and $\delta>\ep$ there exist $m>0$ and $C>0$ and $M$ such that
\begin{equation}
\|\tilde\alpha^{-i\sigma}R(\sigma)\|_{\cL(\alpha^{-1}
\tilde\alpha^\delta H^m_0(\bar X_{1/2}),C^k(\bar X))}
\leq C|\sigma|^M, \label{estimate}
\end{equation}
for $|\sigma|>1$, $\im\sigma<\ep$.
\end{prop}

\section{Asymptotics for solutions of the wave
  equation}\label{section:asymptotics} We now proceed to study the asymptotics of
solutions of the wave equation at $\tfp.$

Suppose $u$ is a solution of the wave equation, $\Box u=0,$ and $u$ is
smooth on $\tilde M^\circ.$ Energy estimates show that $u$ is necessarily
tempered, in the sense that $u\in \rho^{-s}L^2(\bar M)$ for some $s>0.$
This is shown directly below.

Let $\phi=\phi_0(\rho)\in\CI(\bar M)$ be a cutoff function,
with $\rho\in\CIc([0,1))$, identically $1$ near $0.$ If $u$ satisfies $\Box
  u=0$ and $u$ is smooth in $\mu$ (i.e.\ across the side faces), then $\phi
  u$ is smooth in $\mu,$ and
$$
\Box (\phi u)=[\Box,\phi]u=f
$$
where $f$ is also smooth in $\mu,$ and vanishes in a neighborhood of the
temporal face -- since $[\Box,\phi]\in\Diffb^1(\tilde M),$ by
Lemma~\ref{lemma:Box-form}, is supported away from the temporal face.
Moreover, $v=\phi u$ is the {\em unique} solution of $\Box v=f$ in
$\bar M^\circ$ with $v=0$ for $\rho$ sufficiently large.

Now we wish to take the Mellin transform in $T=e^{-t},$ for functions
supported in a neighborhood of the temporal face, namely in $U=\{\rho<1\}.$
Such a neighborhood is equipped with a fibration $\bar M\to X,$ extending
the fibration $M\to X$ in the interior, and there is a natural density
$|dt|=\frac{|dT|}{T}$ on the fibers. In coordinates $(\rhot,\alpha,\omega)
=(T^{\lambda_{\bh}}/\alpha,\alpha,\omega)$ valid near the temporal face
boundary at $r=r_{\bh}$ this density takes the form
$\frac{|d\rhot|}{\lambda_{\bh}\rhot}.$ So, the Mellin transform can be
taken with respect to this fibration and density. Thus the map
$v\mapsto\hat v$ from functions supported near the temporal 
face to functions on $\Omega\times X$, $\Omega\subset\Cx$,
\begin{equation*}
\hat v(\sigma,z)=\int T^{i\sigma} v(T,z)\,\frac{|dT|}{T}.
\end{equation*}
If $v$ is polynomially bounded in $T,$ supported in $T\geq 0$,
with values in a function space $\cH$ in $z$,
this transform gives an analytic
function in a lower half plane (depending on the order of growth of $v$)
with values in $\cH$.

In fact, writing the integral in coordinates valid near the boundary
of the temporal face, $\rhot=T^{\lambda_{\bh}}/\alpha,$
\begin{equation*}
\hat v(\sigma,\alpha,\omega)=\alpha^{i\sigma/\lambda_{\bh}}\lambda_{\bh}^{-1}
\int \rhot^{i\sigma/\lambda_{\bh}}
v(\rhot,\alpha,\omega)\,\frac{|d\rhot|}{\rhot}.
\end{equation*}
The integral is then the Mellin transform of $v$ with respect to $\rhot$
evaluated at $\sigma/\lambda_{\bh}.$ Thus, if $v$ is smooth on $\tilde
M^\circ$, supported in $\{\rhot\in I\}$, $I\subset(0,1)$ compact,
i.e.\ near but not at the temporal face,
then $\hat v$ is in fact analytic in $\Cx$
with values in functions of the form $\tilde\alpha^{i\sigma}\CI$, with
$\CI$ seminorms all bounded by $C_k\langle\sigma\rangle^{-k}$, $k$ arbitrary.
If  $v$ is just supported in $\rhot<1$ and is conormal on $\tilde M,$ then
$\hat v$ is analytic in a lower half plane with values in functions of the
form $\tilde\alpha^{i\sigma}\CI.$

Assuming for the moment that $\phi u$ is polynomially bounded in $T$,
$\Box (\phi u)=f$ becomes $\hat N_b(\Box)\widehat{\phi u}=\hat f$,
where 
\begin{equation*}
\hat N_b(\Box)=
\alpha^{-2}(\sigma^2-\Delta_X)=
\alpha^{-2}\left(\sigma^2-\alpha^2r^{-2}D_r\alpha^2 r^{2}D_r
-\alpha^2r^{-2}\Delta_\omega\right),
\end{equation*}
so
\begin{equation*}
\left(\sigma^2-\Delta_X \right)\widehat{\phi u}=\alpha^2\hat f.
\end{equation*}

If $\phi u$ is polynomially bounded in $T,$ then both $\hat f$ and
$\widehat{\phi u}$ are analytic in $\im\sigma<-C$, and as $f$ is
compactly supported in $\rhot,$ $\hat f$
is in entire analytic with values in functions of the form
$\tilde\alpha^{i\sigma}\CI,$ with $\CI$ seminorms all bounded by
$C_k\langle\sigma\rangle ^{-k},$ $k$ arbitrary.
Thus,
$$
\widehat{\phi u}=R(\sigma)(\alpha^2\hat f),\ \im\sigma<-C,
$$
and we recover $\phi u$ by taking the inverse Mellin transform.

We now return to arbitrary (not a priori polynomially bounded) $u$,
$f=\Box(\phi u)$, as above. Thus,
$\hat f$
is analytic in all of $\Cx$, with values in functions of the form
$\tilde\alpha^{i\sigma}\CI$, with
$\CI$ seminorms all bounded by $C_k\langle\sigma\rangle
^{-k}$, $k$ arbitrary. Now note that the inclusion
$$
\alpha^{1+s}L^\infty(X)\hookrightarrow L^2_0(\bar X_{1/2})
$$
is continuous for every $s>0$. Thus,
the inclusion
$$
\tilde\alpha^{i\sigma-\delta}\alpha^3\CI(\bar X_{1/2})
\hookrightarrow H^m_0(\bar X_{1/2})
$$
is continuous if
$$
\lambda(\ep+\delta)<2,\ \lambda=\max(\lambda_{\bh},\lambda_{\ds}),
$$
and $\im\sigma<\ep$, which is to say
$$
\tilde\alpha^{i\sigma}\alpha^2\CI(\bar X_{1/2})
\hookrightarrow \alpha^{-1}\tilde\alpha^\delta H^m_0(\bar X_{1/2})
$$
is continuous.

In particular, then
$$
\|\alpha^2\hat f\|_{\alpha^{-1}\tilde\alpha^\delta H^m_0(\bar X_{1/2})}
\leq C_k\langle\sigma\rangle^{-k}
$$
for all $k$ in $\im\sigma<\ep<\delta$ (with new constants), $0<\ep<\delta$
sufficiently small.
Proposition \ref{bounds} shows that, for $\ep>0$ sufficiently small
and for all $N$ and $k$,
\begin{equation}\label{eq:res-est-f}
\|\tilde\alpha^{-i\sigma} R(\sigma)(\alpha^2\hat f)
\|_{C^N(\bar X)} \leq C_k|\sigma|^{-k},  \im\sigma<\ep.
\end{equation}
The inverse Mellin transform of $w=R(\sigma)(\alpha^2\hat f)$ is
\begin{equation*}
\check w(T,z)=(2\pi)^{-1}\int T^{-i\sigma}w(\sigma,z)\,d\sigma.
\end{equation*}
Thus,
\begin{equation*}
\check w(\rhot,\alpha,\omega)=(2\pi)^{-1}\int \rhot^{-i\sigma/\lambda_{\bh}}
\alpha^{-i\sigma/\lambda_{\bh}}w(\sigma,\alpha,\omega)\,d\sigma.
\end{equation*}
In view of \eqref{eq:res-est-f} (in particular the analyticity of
$\check w$ in the lower half plane with the stated estimates),
$w=0$ for $T<0$ (as can be seen directly by shifting the contour to
$\im\sigma=-C$, using the off spectrum resolvent estimate
$\|R(\sigma)\|_{\cL(L^2(X,|dh|))}\leq |\im(\sigma^2)|^{-1}$
and letting $C\to +\infty$). Since the unique solution
of $\Box v=f$, $v$ supported in $T\geq 0$, is $\phi u,$ it follows
that $\check w=\phi u.$

Shifting the contour for the inverse Mellin transform for $w$
to $\im\sigma=\ep$
gives a residue term
at $0$, and shows that
\begin{equation*}
\rhot^{-\ep}( \phi u-v)\in L^2([0,\delta)_{\rhot};\CI(\bar X)),
\end{equation*}
where $v$ arises from the residue at $0$, hence is a constant function.
Note that (at the cost of changing $\ep$) this is equivalent
to the analogous statement with $\rho$ replaced by $\rhot$.
The derivatives with respect $t$ satisfy similar estimates. Hence, the
same estimates hold for the conormal derivatives with respect to $\rhot$
(or equivalently $\rho$). We thus deduce the leading part of the
asymptotics of $u$ at the future temporal face, $\tfp.$

This completes the proof of Theorem~\ref{thm:asymptotics}. The main result,
Theorem~\ref{thm:main}, follows from the combination of this result with
Proposition~\ref{prop:interior-propagation}.

\def\cprime{$'$} \def\cprime{$'$}


\end{document}